%% file: main.tex
\documentclass[conference]{IEEEtran}
\IEEEoverridecommandlockouts
\usepackage{cite}
\usepackage{amsmath,amssymb,amsfonts}
\usepackage{graphicx}
\usepackage{textcomp}
\usepackage{xcolor}
\usepackage{color}
\usepackage{listings}
\usepackage[11pt]{moresize}

\usepackage{algorithm}
\usepackage{algpseudocode}
\usepackage{float}
\usepackage{subcaption}
\usepackage{url}
\usepackage{enumerate}

\usepackage{multirow}

\usepackage{todonotes}

\newcommand{\ignore}[1]{}
\newcommand{\ichi}[1]{{\color{brown}#1}}

\begin{document}

\title{An Experimental Study of Two-Level Schwarz Domain Decomposition Preconditioners on GPUs}

\author{
\IEEEauthorblockN{
Ichitaro Yamazaki\IEEEauthorrefmark{2},
Alexander Heinlein\IEEEauthorrefmark{1},
Sivasankaran Rajamanickam\IEEEauthorrefmark{2}
}

\IEEEauthorblockA{
\IEEEauthorrefmark{2}Sandia National Laboratories, Albuquerque, New Mexico, U.S.A\\
\IEEEauthorrefmark{1}Delft University of Technology, Delft, Netherlands}
}
\maketitle

\begin{abstract}
\input{sections/abst}
\end{abstract}

\section{Introduction}
\input{sections/intro}

\section{Related Work}
\input{sections/related}

\section{The GDSW Algorithm}\label{sec:algo}

\input{sections/gdsw}

\section{FROSch software}\label{sec:frosch}
\input{sections/frosch}

\section{GPU Acceleration}
\input{sections/gpu}

\section{Discussion}\label{sec:discussion}
\input{sections/discussion}

\section{Experimental Setup}\label{sec:setups}

\input{sections/setups}

\section{Performance Results}\label{sec:results}
\input{sections/perform}

\section{Conclusion}
\input{sections/concl}

\section*{Acknowledgment}

The authors would like to thank Christian Glusa for helping to integrate {\tt HalfPrecisionOperator} into {\tt FROSch}.
This work was supported in part by the U.S. Department of Energy, Office of Science, Office of Advanced Scientific Computing Research, Scientific Discovery through Advanced Computing (SciDAC) Program through the FASTMath Institute under Contract No. DE-AC02-05CH11231 at Sandia National Laboratories.
Sandia National Laboratories is a multimission laboratory managed and operated by National Technology and Engineering Solutions of Sandia, LLC, a wholly owned subsidiary of Honeywell International, Inc., for the U.S. Department of Energy's National Nuclear Security Administration under contract DE-NA-0003525. This paper describes objective technical results and analysis. Any subjective views or opinions that might be expressed in the paper do not necessarily represent the views of the U.S. Department of Energy or the United States Government.

\bibliographystyle{plain}
\bibliography{ref}
\end{document}

%% file: sections/abst.tex
The generalized Dryja--Smith--Widlund (GDSW) preconditioner is
a two-level overlapping Schwarz domain decomposition (DD) preconditioner
that couples a classical one-level overlapping Schwarz preconditioner
with an energy-minimizing coarse space.
When used to accelerate the convergence rate of Krylov subspace iterative methods, the GDSW preconditioner provides robustness and scalability for the solution of sparse linear systems arising from the discretization of a wide range of partial different equations.
%
%
In this paper, we present {\tt FROSch} (Fast and Robust Schwarz), a domain decomposition solver package which implements GDSW-type preconditioners for both CPU and GPU clusters.
%
%
To improve the solver performance on GPUs,
we use a novel decomposition to run multiple MPI processes on each GPU, 
reducing both solver's computational and storage costs and potentially improving the convergence rate. This allowed us to obtain competitive or faster performance using GPUs compared to using CPUs alone.
We demonstrate the performance of {\tt FROSch} on the Summit supercomputer with NVIDIA V100 GPUs, where we used NVIDIA Multi-Process Service (MPS) to implement our decomposition strategy.

%
The solver has a wide variety of algorithmic and implementation choices, which poses both opportunities and challenges for its GPU implementation. 
We conduct a thorough experimental study with different solver options including the exact or inexact solution of the local overlapping subdomain problems on a GPU.
We also discuss the effect of using the iterative variant of the incomplete LU factorization and sparse-triangular solve as the approximate local solver, and using lower precision for computing the whole {\tt FROSch} preconditioner. 
Overall, the solve time was reduced by factors of about $2\times$ using GPUs, while the GPU acceleration of the numerical setup time depend on the solver options and the local matrix sizes. 

%% file: sections/intro.tex
Domain decomposition methods (DDMs)~\cite{Smith96,Dolean15} may be used to build a class of effective parallel solvers for sparse linear  systems arising from the discretization of partial differential equations. In DDMs, the global problem is decomposed into smaller local subproblems, 
which can be processed in parallel. As a result, DDM preconditioners are well-suited for solving large-scale linear systems on distributed-memory computers. However, for one-level DDMs, the number of iterations required for the solution convergence typically increases with an increasing number of subdomains. As a remedy, a second-level coarse system, which is determined by carefully-designed coarse basis functions, is introduced. As a results  the condition number of the preconditioned matrix, and thus the required number of iterations, becomes asymptotically independent of the number of subdomains.
In this paper, we consider the generalized Dryja--Smith--Widlund (GDSW)~\cite{Dohrmann:2008} two-level Schwarz DDM, which combines classical one-level overlapping additive Schwarz preconditioner with energy-minimizing coarse basis functions and has been shown to be robust and scalable for solving many challenging problems. 

Our focus is on the GPU performance of the GDSW preconditioner.
Several algorithmic and software options are possible for the GDSW algorithm.
Each of these options has multiple tunable parameters, and a good choice of the parameters can be architecture and problem specific.
Some of these options also have algorithmic implications in addition to the implementation choices. For example, 
the computational complexity of the local sparse solver increases more than linearly to the local matrix size.
Since there are typically more CPU cores than the GPUs, decomposing the global domain to one subdomain per GPU instead of one subdomain per CPU core may increase the complexity of the local subdomain solver. Moreover, a fewer subdomains lead to a smaller coarse space, which could degrade the convergence behavior.
Finally, many of the kernels, which DDM solvers depend on, such as sparse direct solver, incomplete factorization, sparse triangular solver, and sparse matrix-matrix multiply, are difficult to optimize on GPUs. 
All of these properties pose challenges when implementing the solver and tuning its performance for the GPU architectures.
As a result, the comprehensive GPU performance study of two-level DDM solver has been lacking, especially at scale, to the best of our knowledge.

\begin{table*}[]
    \centering\footnotesize
    \begin{tabular}{|l|lll|}
    \hline
                            & used                & effects & other options available\\
    \hline\hline
    Krylov solver           & single-reduce GMRES & improve data access      & standard, pipelined, comm-avoid\\
    GDSW                    & two-level rGDSW     & reduce coarse space size & standard, multi-level\\
    Direct solver           & GPU-enabled {\tt Tacho}        & allows use of GPU           & CPU-only, e.g., {\tt SuperLU} \\
    Sparse triangular solve & supernode-based {\tt Kokkos-Kernels}    & improve GPU utilization  & element-based, partitioned inverse\\
    Inexact solver          & iterative variants {\tt FastILU}/{\tt FastSpTRSV} & expose more parallelism  & standard on CPU/GPU\\
    Precision               & single precision {\tt HalfPrecisionOp} & reduce data volume & uniform double precision\\
    \hline
    \# of subdomains        & \# of CPU cores     & reduce solver costs & \# of GPUs\\
                            &                     & \& improve convergence & \\
                            \hline
    \end{tabular}
    \caption{Solver options used in this paper: each solver option has multiple tunable parameters, e.g.,
             {\tt GMRES} with restart length, orthogonalization scheme, {\tt rGDSW} with size of overlap, {\tt Tacho} with matrix ordering scheme, {\tt FastILU} with number of levels, and Jacobi iteration count and damping factor, just to name a few main parameters.}\label{tab:params}
\end{table*}

%
To fill this gap, we study the GPU performance of {\tt FROSch} (Fast and Robust Schwarz), a solver package,
which implements GDSW preconditioners within Trilinos software framework; cf.~\cite{Heinlein:2020}.
Our implementation is also portable to different hardware architectures with a single code base. 
As architectures change rapidly,
it is critical to design the software stack such that the solver is portable
to different hardware architectures (e.g., isolating the hardware-specific codes
and optimizations from the high-level solver design and implementation).
Though this avoids the need of re-writing the solver for each new architecture,
in order to obtain high performance on a specific architecture, including on a GPU cluster,
the software stack must be carefully designed and new variants of the algorithms may be needed.
We evaluate many of the algorithmic and software choices that are critical for the GPU performance of the GDSW preconditioners
including new solver capabilities added for this purpose, 
e.g.,
\begin{itemize}
    \item single-reduce variant~\cite{Katarzyna:2021} of the Krylov solver, which performs only one global-reduce for each iteration,
    \item the exact solution of the local overlapping subdomain and coarse space problems based on the direct sparse LU factorization on a GPU~\cite{Kim:2018}, 
    \item a supernodal-based sparse triangular solver~\cite{Yamazaki:2020}, which reduces the number of kernel launches and exploit the hierarchical parallelism available on a GPU, and
    \item iterative variants of ILU factorization and sparse-triangular solver~\cite{Chow:2015}, which have much higher costs of computation but expose more parallelism than the standard substitution-based algorithm.
\end{itemize}
Table~\ref{tab:params} lists some of the main solver options, which we will evaluate in this paper, and other parameters we selected based on our experience and past results.

Experimental results on the Summit supercomputer with NVIDIA V100 GPUs demonstrate the potential of the two-level DD solvers on the GPU clusters.
We compare the GPU performance with the CPU performance using all the CPU cores on each node. We believe this provides a conservative but fair performance comparison that an application is expected to see.
Using GPUs, the solve time for 3D elasticity problems was reduced by factors of around $2\times$, while
the effects on the numerical solver setup time depends on the solver options and the local matrix sizes. 
In cases of using local direct solvers, the total solution time
(the total of the numerical setup and solve time) for a single linear system was reduced by a factor of about $1.1\times$ to $1.7 \times$. 
If the application requires to solve a sequence of the linear systems with different right-hand-sides,
the cost of the numerical setup can be amortized over multiple solves and the speedups closer to $2\times$ can be obtained.

The main contributions of this work are:
\begin{itemize}
    \item A GPU implementation and large-scale GPU performance study of a two-level DDM solver;
    
    \item A novel decomposition strategy that allows the use of  NVIDIA Multi-Process Service (MPS) to run multiple MPI processes on each GPU, and significantly reduces both the computational and storage cost of the DDM solver, and potentially improves convergence (Section~\ref{sec:discussion}). We are not aware of other studies,
          which use MPS with a production-ready linear solver.
          In our performance studies of using MPS on Summit, both the numerical setup and solve time of {\tt FROSch} was reduced
          by the factors of up to $3\times$.
    \item A detailed experimental study of several solver options for the two-level DDM including direct, incomplete, and approximate factorizations, with multiple parameter choices for each of them.    
    \item Numerical and GPU performance studies with inexact or approximate local linear solver
          or with DDM preconditioner in a lower precision, enabling the solution of larger-scale linear systems 
          than the linear system that the typical DDM solvers can (the exact solution with local direct solvers
          in double precision, which typical GDSW in practice, and its theory, is based on).
\end{itemize}

%% file: sections/related.tex
As GPUs became a critical part in scientific computing,
there have been several works to optimize the computational kernels, which are also needed for DDM solvers.
%
On the other hand, a GPU implementation of a two-level DDM solver that uses these kernels in addition to other kernels has not been demonstrated at scale.
Luo \textit{et al}.~\cite{Luo:2011} investigated the GPU performance of a \emph{one-level} DDM, and the number of MPI processes was limited by the number of GPUs.
We will show that this is sub-optimal in terms of performance and convergence. Solver performance can be improved by using a decomposition that maps one MPI process on each CPU core, and multiple MPI processes on each GPU.
In addition, they used 
 the GPU to only accelerate the local subdomain solver, which was based on smooth-aggregate multi-grid
with \emph{dense} coarse solver. Hence, the GPU performed only sparse-matrix vector multiply, dense vector update, and dense triangular solve,
which are relatively easy to parallelize on a GPU (the paper avoided using ILU since it is difficult to parallelize on a GPU).

In parallel to our work, Šístek and Oberhuber employed GPUs to speed up the local solves in the \emph{two-level} balancing domain decomposition by constraints (BDDC) method in~\cite{sistek_acceleration_2022}; in particular, they perform the factorization and forward and backward solves of \emph{dense} local Schur complement matrices on the GPUs. For time-dependent simulations, where the factorizations can be reused between different time steps, they observed a speedup of up to $5\times$, compared to using CPUs only, by simply storing the dense matrix on GPUs.

In this paper, we will study the GPU performance of the \emph{two-level} DDM solver using the kernels that are more commonly used for the DDM solvers -- \emph{sparse} direct and incomplete factorizations.

%% file: sections/gdsw.tex
We consider the solution of the linear system of equations,
\[
  A x = b.
\]
The DDM solvers have been extensively studied for the matrices
arising from the discretization of an elliptic partial differential equation, in both theory and practice,
but they have been successfully applied to many other problems~\cite{Smith96,Dolean15}.

The two-level overlapping additive Schwarz preconditioner 
is based on a decomposition of the global domain $\Omega$ into $n_p$ nonoverlapping subdomains $\Omega_1, \dots, \Omega_{n_p}$. 
These subdomains are then extended by $\ell$ layers of mesh elements (alternatively mesh nodes) to obtain corresponding overlapping subdomains $\Omega_1', \dots, \Omega_{n_p}'$. 
The preconditioner is then given by
\begin{equation}\label{eq:M}
  M^{-1} = \Phi A_0^{-1} \Phi^T + \sum_{i = 1}^{n_p} R_i^T A_i^{-1} R_i,
\end{equation}
where $R_i$ is the restriction operator from the global domain~$\Omega$ to the $i$th overlapping subdomain $\Omega_i'$ and $A_i = R_i A R_i^T$. 
To construct a robust and efficient preconditioner,
the critical component is the coarse basis functions, the columns of the matrix~$\Phi$, that yields the coarse matrix $A_0 = \Phi A \Phi^T$.

\begin{figure}
\centerline{
\begin{subfigure}[b]{.43\linewidth}
  \includegraphics[width=.88\linewidth]{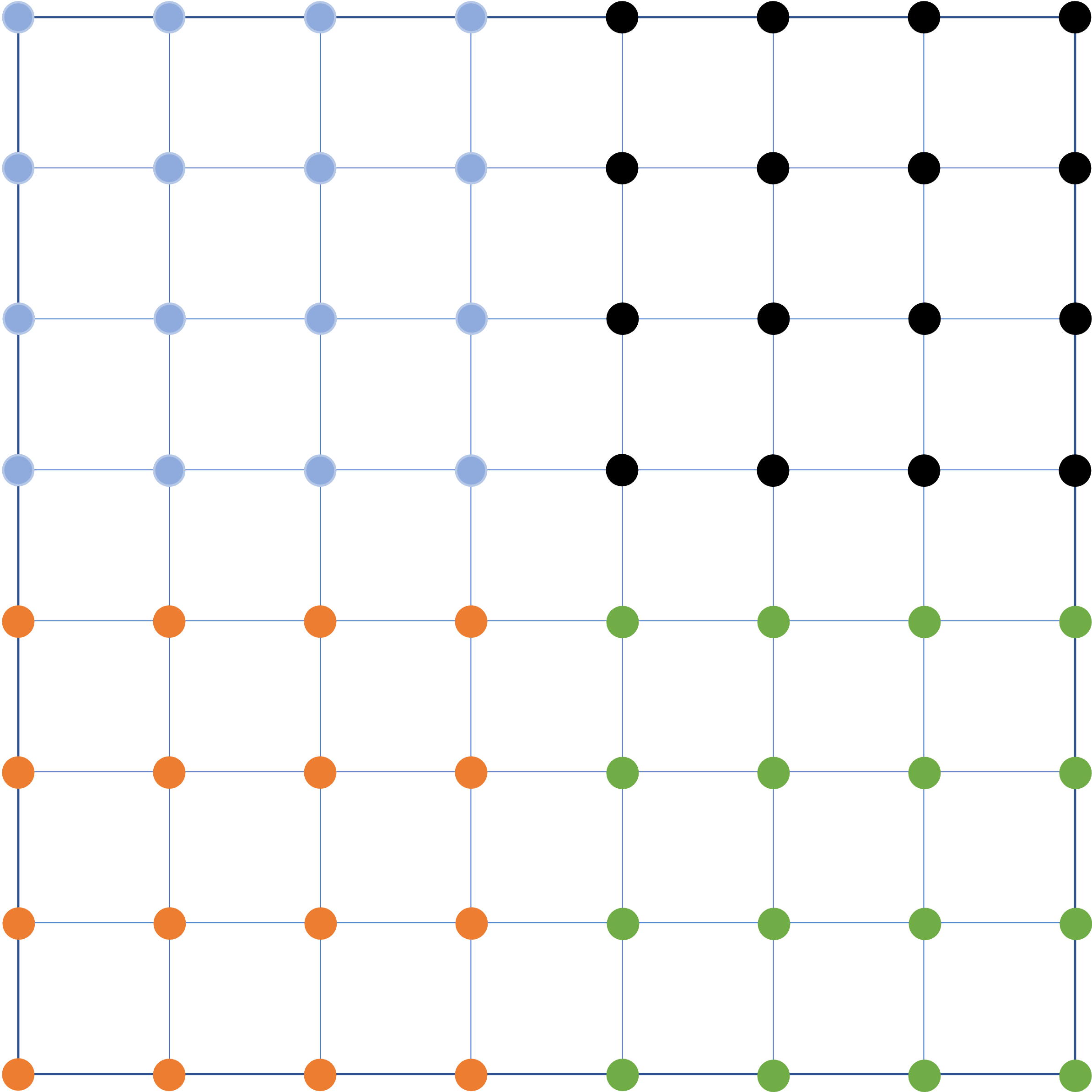}
  \caption{Original problem.}  
\end{subfigure}
\begin{subfigure}[b]{.43\linewidth}
  \includegraphics[width=.9\linewidth]{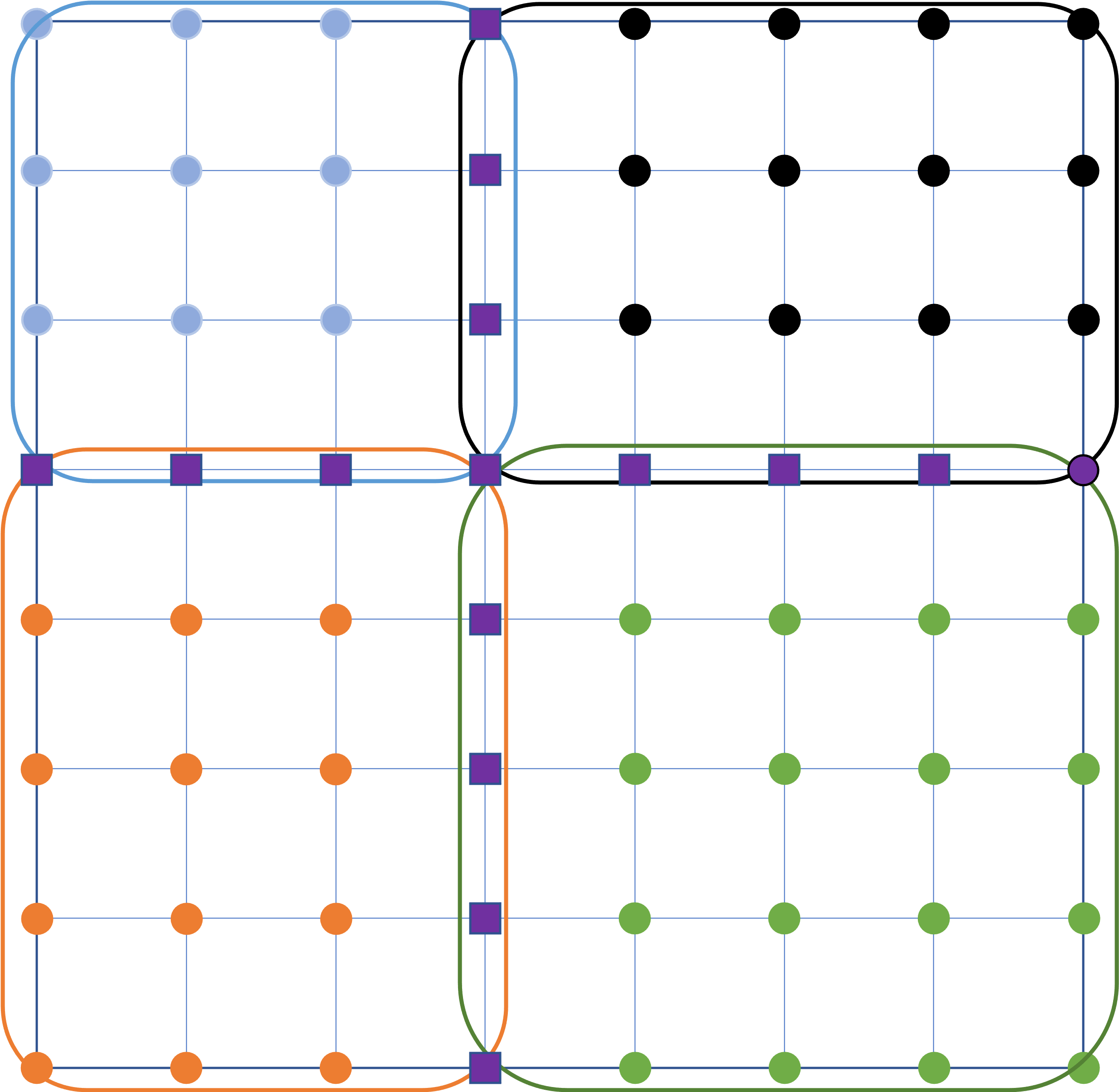}
  \caption{Nonoverlapping partition.}  
\end{subfigure}
}
\caption{Domain decomposition with local $4$-by-$4$ problem.
         In Fig.1(b), interior and interface are represented by circle and square markers,
         respectively.
\label{fig:grid}}
\end{figure}

The coarse basis functions $\Phi$ of GDSW~\cite{Dohrmann:2008} type preconditioners are constructed as energy-minimizing extensions from the interface~$\Gamma$ of the nonoverlapping DD to the interior~$I$ of the subdomains; see Fig.~\ref{fig:grid} for an illustration of the decomposition and see~\cite{Heinlein:2020} for a discussion of the implementation in {\tt FROSch}.
%
For our discussion, we reorder and partition
the global matrix~$A$ into a $2$-by-$2$ block structure
\[
        \left\lbrack \begin{array}{cc}
                    A_{I I} & A_{I \Gamma}\\
                    A_{\Gamma I} & A_{\Gamma \Gamma}
                    \end{array}
        \right\rbrack,
\]
such that the indices $I$ and $\Gamma$ correspond 
interior and the interface degrees of freedom (dofs), respectively.

Let $R_\Gamma$ be the restriction operator
from the global to the interface dofs, such that $A_{\Gamma \Gamma} = R_{\Gamma} A R_{\Gamma}^T$, and 
$n_I$ and $n_\Gamma$ 
denote the numbers of the interior and interface dofs,
respectively.
Then, the GDSW coarse basis functions are defined 
as follows:

1) The interface $\Gamma$ is partitioned into $n_c$ connected components,
$\Gamma_1, \dots, \Gamma_{n_c}$, potentially with overlaps, and $R_{\Gamma_i}$ is the restriction operator from the global interface $\Gamma$ to 
$\Gamma_i$. 

2) 
To obtain a partition of unity on the interface while accounting for the overlapping portions of the interface decomposition, we introduce diagonal scaling matrices~$D_{\Gamma_i}$,
\[
  \sum_{i=1}^{n_c} R_{\Gamma_i}^T D_{\Gamma_i} R_{\Gamma_i} = I_{\Gamma},
\]
where $I_{\Gamma}$ is the identity matrix on $\Gamma$.

3) 
Now, to obtain a robust and efficient preconditioner~$M^{-1}$, the critical component of GDSW type preconditioners is the $n$-by-$n_n$ matrix~$Z$, which contains the null space of the global Neumann matrix corresponding to $A$ as columns.
This matrix may be computed ``algebraically'' for some cases (e.g., just one constant column for a Laplace problem), while in some applications, the null space may be explicitly available.
      
\smallskip
In Section~\ref{sec:results}, we present performance results for a 3D linear elasticity problem, for which, the null space consists of the (linearized) rigid body motions, i.e., translations and linearized rotations. As discussed in~\cite{Heinlein:2021}, the linearized rotations cannot simply be obtained algebraically, however, the method might still perform well when only the translations are used.

4) 
Finally, given the null space matrix $Z$, the energy-minimizing coarse basis functions are computed as
\begin{equation} \label{eq:ext}
	\Phi = \left\lbrack\begin{array}{c}
		-A_{II}^{-1} A_{I\Gamma}\\
		I
	\end{array}\right\rbrack
	\Phi_{\Gamma},
\end{equation}
where $\Phi_{\Gamma}$ is an $n_I$-by-$(n_c n_n)$ matrix given by
\[
  \Phi_{\Gamma} = \lbrack R_{\Gamma_1}^T \Phi_{\Gamma_1}, \;\dots,\;
                          R_{\Gamma_{n_c}}^T \Phi_{\Gamma_{n_c}} \rbrack
\]
and each $n_I$-by-$n_n$ matrix $\Phi_{\Gamma_i}$ spans the null space restricted to the $i$th interface,
$D_{\Gamma_i} R_{\Gamma_i}(R_\Gamma Z)$.
Hence, $\Phi_{\Gamma}$ 
has dimension $n_I$-by-$(n_n n_c)$, while the dimension
of $\Phi$ is $n$-by-$(n_n n_c)$.

The computation of the subdomain problems in~\eqref{eq:M} parallelizes well since, in a distributed-memory implementation with MPI, the $i$th subdomain $A_i$ is assigned to the $i$th MPI process and can be processed in parallel. The extensions~\eqref{eq:ext} 
can be parallelized similarly since $A_{II}$ has a block diagonal structure, $A_{II} = \text{diag} ( A_{I_1 I_1}, \ldots, A_{I_{n_p}, I_{n_p}} )$,
where 
$A_{I_i I_i}$ corresponds to the
interior part of the $i$th nonoverlapping subdomain. 
%
The GDSW 
coarse space can keep the condition
number of the preconditioned matrix $AM^{-1}$, and hence the number
of iterations, asymptotically constant with an increasing number of subdomains; see, e.g.,~\cite{Dohrmann:2008}.

There are several variants of GDSW including:
\begin{itemize}
\item 
As the number of subdomains
and MPI processes increases, the solution of the coarse problem $A_0^{-1}$ eventually becomes a parallel performance bottleneck. 
To alleviate this bottleneck
a ``reduced'' variant of GDSW (rGDSW) only uses only coarse basis functions
corresponding to vertices, but not to faces or edges; see~\cite{Dormann:2017,Heinlein:2020:RDG}. Moreover, multi-level approaches have been proposed to recursively apply GDSW on the coarse problem; cf.~\cite{Heinlein:2021:Copper}. 

\item
To enhance the coarse space for problems with a highly heterogeneous coefficient, potentially with high jumps, ``adaptive'' GDSW (AGDSW) enriches the coarse space by additional components that are computed
by solving local generalized eigenvalue problems; see, e.g.,~\cite{Heinlein:2019}.

\end{itemize}
We do not use the adaptive or multi-level variants in this paper, however several of our results from this study apply to these variants as well. 
In addition, the behavior of GDSW has been extensively compared to other two-level DDMs, and many of the DDM solvers use similar underlying kernels. Hence our experimental study may provide insights to other methods.

%% file: sections/frosch.tex
\begin{figure*}[h]
\centerline{
\begin{subfigure}[b]{.25\linewidth}
  \centering
  \includegraphics[width=\linewidth]{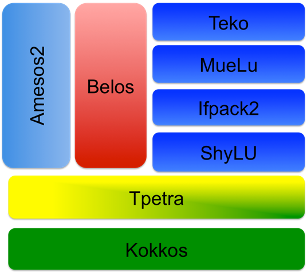}
  \caption{Trilinos solver software stack.}  
  \label{fig:solver-stack}
\end{subfigure}
\quad\quad\quad
\begin{subfigure}[b]{.6\linewidth}\footnotesize
  \begin{tabular}{ll}
  \multicolumn{2}{l}{Flexible Solver Interface} \\
  \hline
  {\tt ShyLU}   & Distributed DD preconditioner ({\tt FROSch}) and\\
                & on-node factorization-based local solvers ({\tt Basker}, {\tt Tacho})\\
  {\tt MuLue}   & Algebraic multigrid solver\\
  {\tt Amesos2} & Direct solver interfaces (e.g., {\tt KLU}, {\tt PaRDISO}, {\tt SuperLU}, {\tt Tacho})\\
  {\tt Belos}   & Krylov solvers (e.g., CG, GMRES, BiCG, and \\
                & their communication-avoiding or pipelined variants)\\
  {\tt Ifpack2} & Algebraic preconditioners (ILU, relaxation, one-level Schwarz) \\
  \\
  \multicolumn{2}{l}{Portable Performance} \\
  \hline
  {\tt Tpetra}         & Distributed sparse/dense matrix-vector operations \\
  {\tt Kokkos-Kernels} & Performance portable on-node graph and sparse/dense matrix operations\\
  {\tt Kokkos}         & C++ programming model for performance portable applications \\ 
                       & on different node architectures (e.g., CPUs, NVIDIA/AMD GPUs)
  \end{tabular}
\caption{Linear solver package descriptions.}
\label{fig:solver-package}
\end{subfigure}
}
\caption{Trilinos linear solver packages}\label{fig:trilinos}
\end{figure*}

{\tt FROSch}~\cite{Heinlein:2020}
implements GDSW type preconditioners within the {\tt Trilinos} software framework~\cite{trilinos-website},
a collection of open-source software packages that can be used as building blocks
for developing large-scale scientific applications.
Fig.~\ref{fig:trilinos} shows the core {\tt Trilinos} packages
for solving linear systems of equations.
These packages can be combined to develop a flexible and adaptable solver for large-scale scientific applications.
For instance, {\tt FROSch} has interfaces to these solver packages
for solving its local overlapping subdomain and coarse problems:
direct solvers ({\tt Amesos2}~\cite{Bavier:2012}), 
inexact and preconditioned Krylov solvers ({\tt Ifpack2}~\cite{Ifpack2} and {\tt Belos}~\cite{Bavier:2012}),
and even a local algebraic multigrid solver
or Schwarz methods ({\tt MueLu}~\cite{MueLu} and {\tt FROSch}~\cite{Heinlein:2020}).
In addition, {\tt FROSch} can be used as a preconditioner for {\tt Belos},
which implements 
Krylov solvers, including variants 
which can be optimized for the GPU architectures,
such as single-reduce, communication-avoiding, and pipelined variants~\cite{Yamazaki:2020:PP}.


Furthermore, {\tt FROSch} builds on the {\tt Trilinos} software stack, specifically packages that
provide the portable performance on different hardware architectures:
In particular, 
\begin{itemize}
    \item {\tt Kokkos}~\cite{Trott:2021} is a C++ performance-portable programming ecosystem.
          It provides the memory abstraction and functionality to dispatch particular functions
          for parallel-operations on a specific execution space 
          on a CPU or GPU. This enables portable thread performance
          on different manycore architectures using a single code base (assuming algorithms are performance portable).
    \item {\tt Kokkos-Kernels}~\cite{KokkosKernels} is a collection of Kokkos-based kernels for on-node sparse or dense matrix, or graph operations on CPUs and GPU.
    \item {\tt Tpetra}~\cite{tpetra-website} implements distributed graph, matrix, and vector operations
          for CPU and GPU clusters.
\end{itemize}
Though we focus on the {\tt FROSch} software stack, this is not the only option
for the portable performance stack.
For instance,
previous studies have compared the performance of the on-node portable layers
and individual kernels~\cite{Trott:2021}.

%% file: sections/gpu.tex
Many of the current high-performance computers
are composed of the heterogeneous compute node architectures,
i.e., each node consists of multicore CPUs and multiple GPU accelerators.
A GPU with a large number of compute cores and a high memory bandwidth
is suited for a highly-parallel computation with a regular memory access pattern,
while some operations are better suited for CPUs. 

This poses both challenges and opportunities for designing high performance 
sparse linear solvers, including {\tt FROSch}, where most of the required computational kernels
have irregular memory accesses and a small ratio of the computation to the data accesses.
As a result, their performance is often bounded not by the computation
but by the memory bandwidth, if not by memory latency.
In order for the solver to utilize the GPUs well, 
both the solver and its underlying computational
kernels must be carefully designed, and new variants of the algorithms may need to be developed. 
In this section, we discuss some of the specific approaches taken to improve the performance of {\tt FROSch}
on GPU clusters.

\subsection{Software Considerations}

\ignore{
\subsubsection{UVM removal}

With a CUDA GPU, UVM simplifies the coding by providing automatic data migration between the CPU and GPU memories.
However, the software needs to be carefully designed, and it is still difficult to manage the data migration
and obtain high performance.
Moreover, there are concerns about how well this will be supported, especially in term of performance, on future GPU systems.
{\tt Trilinos} software stack has been updated to provide code path without UVM
but instead with explicit data migration between the CPU and GPU memories~\cite{Tpetra:UVM:Multivector,Tpetra:UVM:CrsMatrix}.
\ichi{...}}

\subsubsection{Software Structure}

In many scientific and engineering simulations, we often
perform the numerical factorization multiple times
for a given mesh with the same sparsity structure or need 
to solve a sequence of linear systems with different right-hand-side
vectors. 
Hence, all the linear solvers
in {\tt Trilinos} have three distinct phases:
\begin{enumerate}[(a)]
    \item {\it Symbolic Factorization}, given a sparsity structure or a graph, 
          performs all the symbolic analysis and factorization
          and allocate required GPU memories. Operations such as the symbolic analysis for an LU factorization, computing the level sets for a triangular solve, are done here. This is typically done on a CPU.
    \item {\it Numerical Factorization}, given numerical values of the input matrix,
          performs the numerical factorization; in {\tt FROSch}, this part includes the computation of the coarse basis functions, computing the coarse space matrix, and factoring the overlapping local subdomain and coarse matrices. Steps such as the sparse matrix - sparse matrix multiplication for computing $A_0$, numerical factorization of the LU factorization or incomplete factorization are also part of this step. We compute these on the GPUs, when appropriate.
    \item {\it Solve Phase}, given right-hand-side vector(s), compute the solution to the linear system. The sparse triangular solve for the direct or incomplete 
        factorization of local matrix is done on the GPU as part of this phase.
\end{enumerate}
These distinct phases are critical, especially for GPUs
since large parts of the symbolic analysis are difficult to parallelize,
and the GPU memory allocations can take a significant amount of time.

\subsubsection{Lower Precision Preconditioning}

Within {\tt Trilinos}, the software package {\tt Belos} implements Krylov solvers.
It uses the {\tt Operator} class for applying a preconditioner, including
algebraic multigrid (AMG) and domain decomposition (DD) preconditioners,
which are implemented in the {\tt MueLu} and {\tt FROSch} packages, respectively.
These preconditioners are typically constructed from a sparse matrix class,
called {\tt CrsMatrix}.

A new utility function that converts
a {\tt CrsMatrix} object into a new object in half the precision
was developed (e.g., if the original matrix is in double precision, then
the new matrix will be in single precision), allowing users to construct
the preconditioner also in half the precision.
The new {\tt HalfPrecisionOperator} class, 
which inherits the base {\tt Operator} class in the working precision
and internally holds the operator in half the precision as its member variable,
is also implemented. 
When this new operator is applied to vectors,
it internally type-casts the input vectors into half the precision,
applies the operator (e.g., preconditioner) in half the precision, and then type-casts the resulting output vector back
into the working precision. Though it has the overhead of type-casting,
these new capabilities allow users to apply many of the preexisting Trilinos preconditioners
in half the precision within the current Trilinos framework.
Currently, {\tt MueLu} and {\tt FROSch} (whose cost of applying is typically
much higher than that of type-casting vectors) have been extended to utilize
these new capabilities

\subsection{Algorithmic Consideration}

\subsubsection{Sparse Direct Matrix Factorization}

Most of the theoretical results for DD solvers
(including the condition number estimates for the preconditioned matrix)
assume the exact solution of the overlapping subdomain
and coarse problems.
As a result, in practice, DD solver typically use sparse direct solvers.

Sparse direct solvers are a critical component in many scientific applications,
and there have been extensive efforts to develop high performance
sparse direct solvers~\cite{Davis:2016}.
For our experiments, we used {\tt SuperLU} and {\tt Tacho}
software packages that provide two different approaches to the sparse direct solvers:
\begin{itemize}
    \item {\tt SuperLU}~\cite{Demmel:1999} implements left-looking sparse LU factorization with partial pivoting.
          It mainly targets a single CPU core, though
          it could be linked to threaded BLAS or LAPACK for runing on multicore CPUs.
          It uses the supernodal block structures of the LU factors in order to
          exploit the memory hierarchy.
    \item {\tt Tacho}~\cite{Kim:2018} is based on multifrontal factorization with pivoting only inside the frontal matrices.
          The original {\tt Tacho} used the {\it task} programming model of {\tt Kokkos}. Though the current implementation
          still uses {\tt Kokkos},  it exploits the hierarchical parallelism
          available on a GPU through the combination of level-set scheduling and team-level BLAS/LAPACK like kernels. 
          It also has interface to the vendor-optimized kernels (i.e., NVIDIA's CuBLAS/cuSolver and AMD's rocBLAS/rocSolver)
          to factorize large frontal matrices with GPU streams. {\tt Tacho} currently supports Cholesky, LDL$^T$, or LU factorization
          of a symmetric positive definite, symmetric indefinite, or numerically nonsymmetric matrix but with symmetric pattern, respectively.
\end{itemize}

\subsubsection{Sparse Triangular Solve}
When a direct sparse matrix factorization is used, 
the resulting sparse triangular matrix typically has the
dense blocks called {\it supernodes}. 
It is possible to exploit this supernodal structure 
to accelerate the triangular solves.
For instance, {\tt Kokkos-Kernels} implements
sparse-triangular solver based on level-set scheduling
of supernodal blocks~\cite{Yamazaki:2020}.
Working with blocks instead of matrix elements
may give several performance advantages on a GPU.
For instance,
it reduces the height of the level-set trees
and the length of the critical path in the parallel execution
(e.g., number of kernel launches) of the sparse triangular solve.
In addition, it allows the hierarchical parallelization,
which fits well to the hierarchical parallelism available on a GPU and
can be exposed using team-based kernels in {\tt Kokkos-Kernels}.

The sparse-triangular solver in {\tt Kokkos-Kernels}
also has an option to perform the partitioned inverse~\cite{Alvarado:1993}
that transforms the sparse triangular solve into a sequence
of sparse-matrix vector multiply which provides more parallelism
than the standard substitution-based algorithm~\cite{Yamazaki:2020}.

\subsubsection{Incomplete Local Solver}

Though DD theory is based on exact solution of the local subdomain and coarse problems, an inexact local solver may work well in practice, in particular, if its application is somewhat spectrally equivalent to applying an exact solver.
In this paper, we explore inexact local solver based on \emph{level}-based 
incomplete sparse LU factorization.
Though several parallel ILU implementations have been proposed,
the standard paralelization scheme for the ILU and spares-triangular solve
is based on the level-set scheduling~\cite{Anderson:1989}.
In Trilinos, these are implemented 
as {\tt SpILU} and {\tt SpTRSV} in {\tt Kokkos Kernels}~\cite{KokkosKernels}.

Though incomplete factorization leads to a fewer fills
and may expose more parallelism,
it may still not provide enough parallelism to utilize a GPU.
To expose more parallelism,
iterative variants of sparse approximate factorization
and of sparse triangular solver are proposed~\cite{Chow:2015}.
It uses Jacobi iterations to approximate each entries of the LU factors
or to approximately solve the linear system with a sparse triangular matrix.
Though each iteration requires the about same number of floating point
operations as the standard algorithms,
this variant has significantly more parallelism.
As a result, when the solver needs a small number of iterations to obtain the solution
of the desired approximation accuracy (our default is five iterations for both), 
it may obtain much shorter time to solution on a GPU.
In Trilinos, these are implemented as {\tt FastILU} and {\tt FastSpTRSV}~\cite{Boman:2016}.


%% file: sections/discussion.tex
\begin{figure}
\centerline{
\begin{subfigure}[b]{.4\linewidth}
  \includegraphics[width=.9\linewidth]{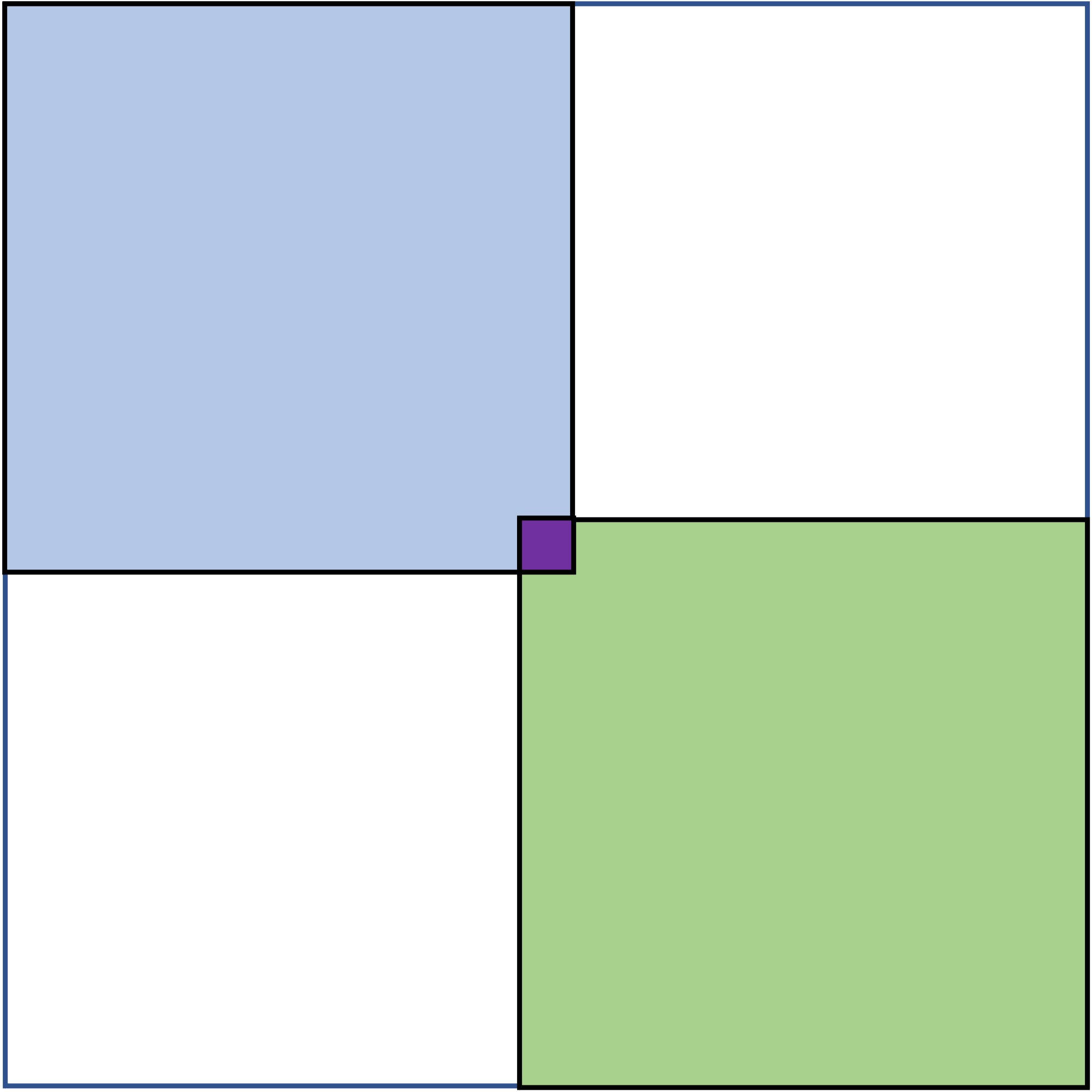}
  \caption{two subdomains.}  
\end{subfigure}
\quad
\begin{subfigure}[b]{.4\linewidth}
  \includegraphics[width=.9\linewidth]{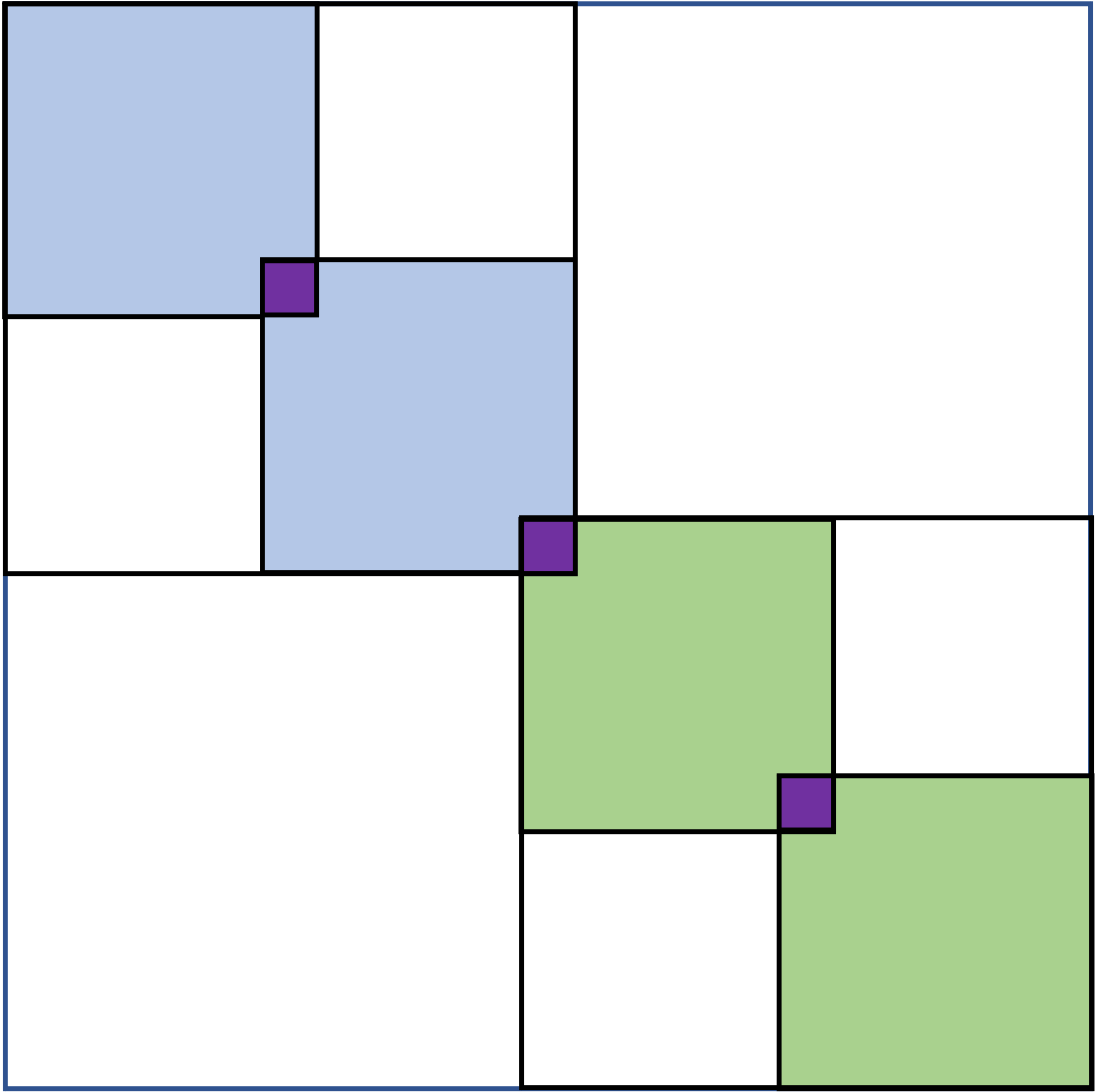}
  \caption{Four subdomains.}  
\end{subfigure}
}
\caption{Overlapping domain decomposition into two or four subdomains (e.g., one or two MPIs on each of two GPUs).
 Instead of assigning two large subdomains to two GPUs (left), we assign four small subdomains to two GPUs (right).\label{fig:dd}}
\end{figure}


Current heterogeneous node architectures typically have more CPU cores than GPUs on each node
(e.g., each node of the Summit supercomputer has 42 IBM Power9 CPU cores and 6 NVIDIA V100 GPUs).
Such a heterogeneous node architecture poses challenges for a node-to-node performance comparison
of the DD solver on CPUs and with GPUs.
This is because, in most cases,
the computational complexity of the local sparse solver increases more than linearly to the local matrix size.
For example, for a 3D problem, when the nested dissection~\cite{George:1973} is used to permute the local matrix with dimension of $n_i$, the sparse direct factorization and corresponding sparse-triangular solve of the local problem typically have the computational complexities of $\mathcal{O}(n_i^2)$ and $\mathcal{O}(n_i^{4/3})$, respectively. As a result, for our strong parallel-scale studies, the computational cost of the DD solver decreases superlinearly with the number of MPI processes. Since each node typically has fewer GPUs than CPU cores, if we launch one MPI process on each CPU core for CPU runs and one MPI process on each GPU for GPU runs (these are the most common setups in practice), each process has a much smaller computational cost for the CPU runs than for the GPU runs.

Moreover, the condition number of the matrix preconditioned with a GDSW preconditioner is bounded as follows:
\begin{align*}
    \kappa \left( M_{\rm GDSW}^{-1} A \right) \leq C
    \left(1+\frac{H}{\delta}\right) \left( 1 + \log\left(\frac{H}{h}\right)\right)^2;
\end{align*}
where $H$ is the maximum diameter of the subdomains and $\delta$ is the width of the overlap
cf.~\cite{Dohrmann:2008}. Hence, the condition number will decrease with a smaller subdomain size ($H/h$), e.g.,
as the number of subdomains increases with a fixed problem size. 

For our experiments, we used NVIDIA Multi-Process Service (MPS) to run multiple MPI processes on each GPU
(see Fig.~\ref{fig:dd}).
Compared to having just one process on each GPU,
this not only reduces the computational and storage costs of the DD solver,
but it may also improve the condition number of the preconditioned matrix, and hence the convergence rate of the Krylov solver.
It is possible to obtain the same decomposition by
having multiple subdomains per MPI and using GPU streams. 
However, this will require significant algorithmic innovations and software efforts for a two-level solver. 
Though it may not be optimal, 
MPS allows us to run the existing code without these code changes
and can provide significant performance gain running multiple subdomains per GPU.

%% file: sections/setups.tex
For all of our performance results presented in this paper,
we used the ``reduced'' GDSW coarse space with an algebraic overlap of one .
Then, as our Krylov solver,
the single-reduce variant~\cite{Katarzyna:2021} of the Generalized Minimum Residual (GMRES) method~\cite{Saad:1986} was used, 
which is a popular Krylov method for solving nonsymmetric linear systems of equations.
We used the restart length of 30,
and considered GMRES to be converged when the residual norm is reduced by a factor of $10^{-7}$.
Finally, we focused on solving 3D elasticity problems in this paper.
Though there are several other preconditioning options in {\tt FROSch},
since our focus is on the performance comparison, and not on the numerical study of GDSW,
these setups provide representative performance of {\tt FROSch}.

We present performance results on the Summit Supercomputer at Oak Ridge Leadership Computing Facility.
Each node of Summit has 42 IBM Power9 CPU cores and 6 NVIDIA V100 GPUs.
Unless specified otherwise, for our CPU runs, we launched 42 MPI processes on each node (one MPI per CPU core), while
for our GPU runs, we used NVIDIA Multi-Process Service (MPS) to run up to 7 MPI processes on each GPU (up to 42 MPI processes per node).
The codes were compiled using {\tt CUDA 10.2.89} and {\tt GCC 7.5.0},
and linked to the vendor-optimized libraries, NVIDIA's {\tt CUBLAS}, {\tt CuSparse} on GPUs, and 
IBM's Engineering and Scientific Subroutine Library {\tt (ESSL) 6.3}.

%% file: sections/perform.tex
\begin{table}[t]
\centerline{
\begin{subfigure}[b]{\linewidth} \centering \scriptsize \addtolength{\tabcolsep}{-0.2mm}
	\begin{tabular}{|ll|rrrrr|}
		\hline
		\multicolumn{2}{|l|}{\# comp.~nodes}                               &                    1 &                    2 &                    4 &                    8 &                   16 \\ \hline\hline
		\multicolumn{2}{|l|}{matrix size}                                  &                 375K &                 750K &                 1.5M &                   3M &                   6M \\ \hline
		\multicolumn{2}{|l|}{CPU}                                          &      {\bf 2.03 (75)} &      {\bf 2.07 (69)} &      {\bf 1.87 (61)} &      {\bf 1.95 (58)} &      {\bf 2.48 (69)} \\ \hline
		\multirow{4}{*}{\rotatebox[origin=c]{90}{GPU}} & $n_p$/gpu = 
   1 &            1.43 (47) &            1.52 (53) &            2.82 (77) &            2.44 (68) &            2.61 (75) \\
		                                               & 2                 &            1.03 (46) &            1.36 (65) &            1.37 (60) &            1.52 (65) &            1.98 (86) \\
		                                               & 4                 &            0.93 (59) &            0.91 (53) &            0.98 (59) &            1.33 (77) &            1.21 (66) \\
		                                               & 6                 &            0.67 (46) &            0.99 (65) &            0.92 (57) &            0.91 (57) &            0.95 (57) \\
		                                               & 7                 &      {\bf 1.03 (75)} &      {\bf 1.04 (69)} &      {\bf 0.90 (61)} &      {\bf 0.97 (58)} &      {\bf 1.18 (69)} \\ \hline
		\multicolumn{2}{|l|}{\bf speedup}                                  & $\mathbf{2.0\times}$ & $\mathbf{2.0\times}$ & $\mathbf{2.1\times}$ & $\mathbf{2.0\times}$ & $\mathbf{2.1\times}$ \\ \hline
	\end{tabular}
   \caption{SuperLU.\\\;\smallskip}
\end{subfigure}}
\centerline{
\begin{subfigure}[b]{\linewidth} \centering \scriptsize \addtolength{\tabcolsep}{-0.2mm}
   \begin{tabular}{|ll|rrrrr|}
   	\hline
   	\multicolumn{2}{|l|}{\# comp.~nodes}                               &                    1 &                    2 &                    4 &                    8 &                   16 \\ \hline\hline
   	\multicolumn{2}{|l|}{matrix size}                                  &                 375K &                 750K &                 1.5M &                   3M &                   6M \\ \hline
   	\multicolumn{2}{|l|}{CPU}                                          &      {\bf 1.60 (75)} &      {\bf 1.63 (69)} &      {\bf 1.49 (61)} &      {\bf 1.51 (58)} &      {\bf 1.90 (69)} \\ \hline
   	\multirow{4}{*}{\rotatebox[origin=c]{90}{GPU}} & $n_p$/gpu = 
   1 &            1.17 (47) &            1.37 (53) &            1.92 (77) &            1.78 (68) &            2.21 (75) \\
   	                                               & 2                 &            0.79 (46) &            1.14 (65) &            1.05 (60) &            1.18 (65) &            1.70 (86) \\
   	                                               & 4                 &            0.85 (59) &            0.81 (53) &            0.78 (59) &            1.22 (77) &            1.19 (66) \\
   	                                               & 6                 &            0.60 (46) &            0.86 (65) &            0.75 (57) &            0.84 (57) &            0.91 (57) \\
   	                                               & 7                 &      {\bf 0.99 (75)} &      {\bf 0.93 (69)} &      {\bf 0.82 (61)} &      {\bf 0.93 (58)} &      {\bf 1.22 (69)} \\ \hline
   	\multicolumn{2}{|l|}{\bf speedup}                                  & $\mathbf{1.6\times}$ & $\mathbf{1.8\times}$ & $\mathbf{1.8\times}$ & $\mathbf{1.6\times}$ & $\mathbf{1.6\times}$ \\ \hline
   \end{tabular}
   \caption{Tacho.}
\end{subfigure}
}
\caption{Total iteration time in seconds and Iteration count for Weak-scale 3D elasticity problems on Summit. \label{tab:solve}}
\end{table}

\subsection{Exact Local Solvers}\label{sec:results-exact}

In this section, we study the performance of {\tt FROSch}
using exact solution of the local overlapping subdomain and coarse space problems.
The nested dissection ordering from {\tt Metis}~\cite{Karypis:2013} was used to
reduce the number of fills in the LU factors, and also to expose more parallelism.
We used either {\tt SuperLU} or {\tt Tacho} to factor our local and coarse matrices on CPUs or GPUs, respectively.
To apply the preconditioner on a GPU with the LU factors computed by {\tt SuperLU}, the supernodal sparse-triangular solver~\cite{Yamazaki:2020} from {\tt Kokkos-Kernels} was used, 
while on CPU, we used the {\tt SuperLU's} internal triangular solver since {\tt Kokkos-Kernels} solver is designed to exploit
the manycore architectures and is not suited on a single CPU core.
With {\tt Tacho}, we used its internal sparse-triangular solver for both CPU and GPU runs.
We did not use the partitioned inverse, and all the sparse-triangular solvers, either on a CPU or on a GPU, 
are numerically equivalent.

Table~\ref{tab:solve} shows the the weak-parallel scaling of the total iteration time required
for the solution convergence, where the local problem size on each compute node is fixed and the global matrix
size grows linearly to the number of the compute nodes.
We used MPS to run multiple MPI processes on each GPU.
As we discussed in Section~\ref{sec:setups}, with the direct factorization of of the local overlapping subdomain matrix,
the resulting sparse-triangular solve has the computational cost that scales superlinearly to the size of the local matrix. 
\emph{Hence, as we map more MPI processes on each GPU using MPS, the local subdomain becomes smaller, 
and the iterative solution time is reduced, significantly (with speedups of $1.3 \sim 2.7\times$).}
Overall, using GPUs, the solution time was reduced by a factor of around $2\times$, compared to the CPU runs.

In some applications, the setup time can also take a significant share of the total simulation time.
Hence, we now study the numerical setup time of {\tt FROSch}.
Fig.~\ref{fig:breakdown} shows the breakdown of the numerical setup time on a single compute node of Summit.
As we expect, especially on CPUs, a significant part of the numerical setup time is spent by the sparse direct solver.
For both CPU and GPU runs with {\tt SuperLU}, the local overlapping and coarse matrices are factored on CPU,
and the factorization time are the same on CPUs and with GPUs. \emph{On the other hand, {\tt Tacho} can exploit the GPU,
and the local factorization time was reduced for the GPU run by $2.4\times$}. This is the first benefit of using GPUs.
Unfortunately, we also see that
some of the setup time beside the sparse direct solver is running slower with GPUs (``black'' part of the bar)\footnote{%
This is mostly due to sparse-sparse matrix product to form the coarse matrix and communication to form the local overlapping subdomain matrix.},
and a significant amount of time is spent setting up the {\tt Kokkos-Kernels} sparse-triangular solve with {\tt SuperLU}: 

\begin{figure}[t]
 \centerline{
   \begin{subfigure}[b]{.48\linewidth}
   \centering
   \includegraphics[width=\linewidth]{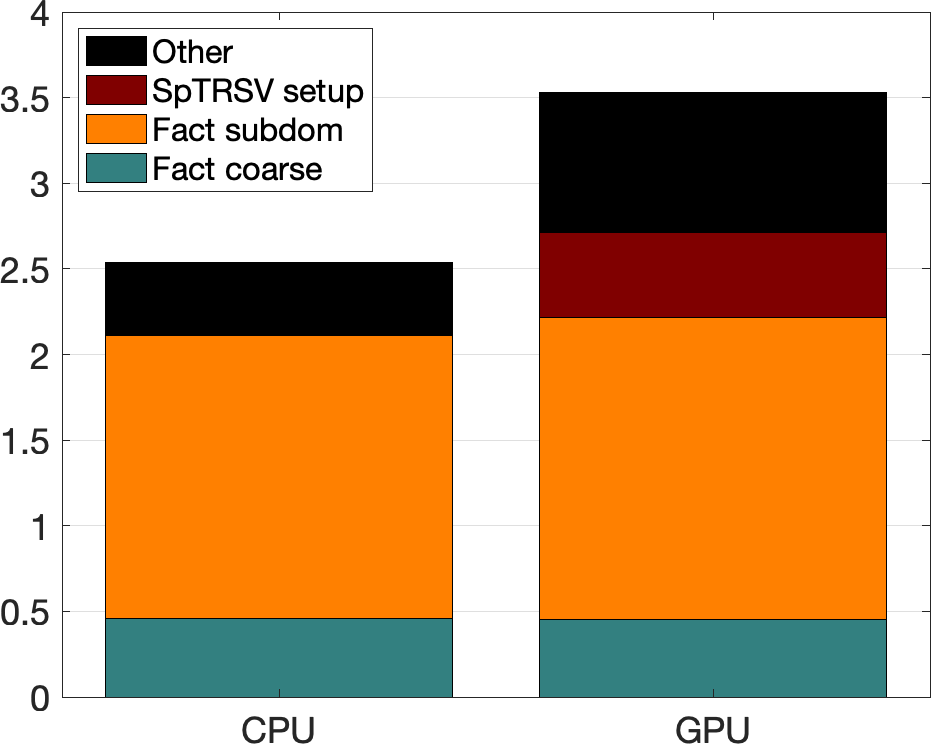}
   \caption{SuperLU}
   \end{subfigure}
 \quad
   \begin{subfigure}[b]{.48\linewidth}
   \centering
   \includegraphics[width=\linewidth]{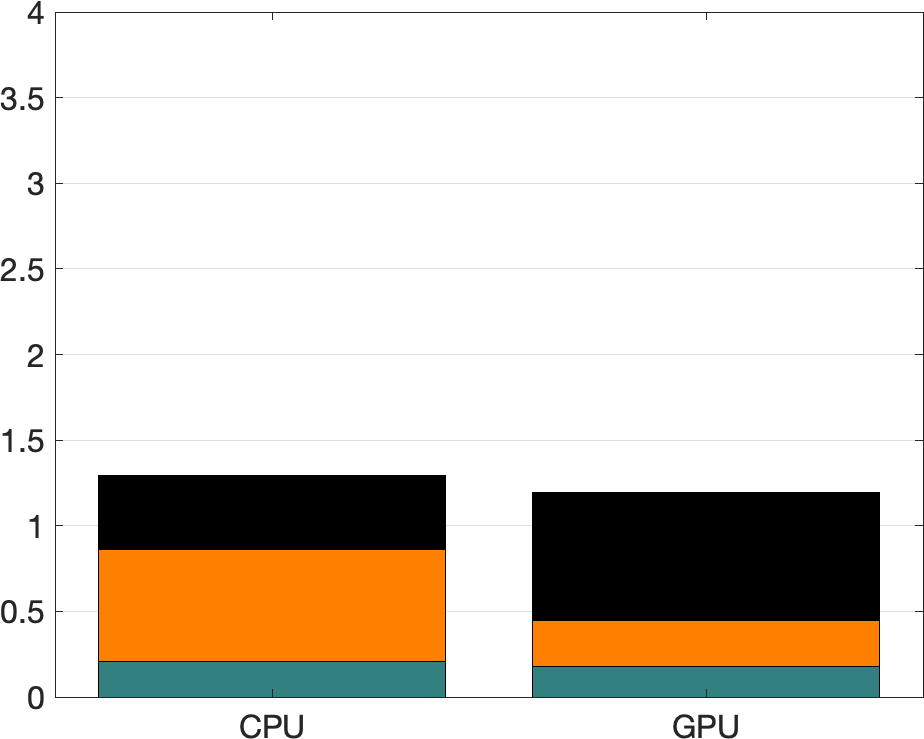}
   \caption{Tacho}
   \end{subfigure}
 }
 \caption{Breakdown of the numerical setup time on one node of Summit ($n=375$K on 42 MPI processes).
 }
 \label{fig:breakdown}
\end{figure}

\begin{table}[t]
\centerline{
\begin{subfigure}[b]{\linewidth} \centering \footnotesize 
   \begin{tabular}{|ll|rrrrr|}
   	\hline
   	\multicolumn{2}{|l|}{\# comp.~nodes}                              &                    1 &                    2 &                    4 &                    8 &                   16 \\ \hline\hline
   	\multicolumn{2}{|l|}{matrix size}                                 &                 375K &                 750K &                 1.5M &                   3M &                   6M \\ \hline
   	\multicolumn{2}{|l|}{CPU}                                         &            {\bf 2.5} &            {\bf 3.0} &            {\bf 3.3} &            {\bf 3.8} &            {\bf 3.6} \\ \hline
   	\multirow{4}{*}{\rotatebox[origin=c]{90}{GPU}} & $n_p$/gpu = 
  1 &                 60.0 &                 71.5 &                 70.9 &                 85.6 &                 85.1 \\
   	                                               & 2                &                 22.9 &                 22.3 &                 26.1 &                 26.9 &                 25.8 \\
   	                                               & 4                &                  8.4 &                  9.4 &                  9.1 &                  9.8 &                 10.2 \\
   	                                               & 6                &                  5.5 &                  5.2 &                  5.2 &                  6.2 &                  6.3 \\
   	                                               & 7                &            {\bf 3.5} &            {\bf 4.2} &            {\bf 4.8} &            {\bf 5.4} &            {\bf 5.4} \\ \hline
   	\multicolumn{2}{|l|}{slowdown}                                    & $\mathbf{1.4\times}$ & $\mathbf{1.4\times}$ & $\mathbf{1.5\times}$ & $\mathbf{1.4\times}$ & $\mathbf{1.5\times}$ \\ \hline
   \end{tabular}
   \caption{SuperLU.\\\;\smallskip}
\end{subfigure}
}
\centerline{
\begin{subfigure}[b]{\linewidth} \centering \footnotesize 
	\begin{tabular}{|ll|rrrrr|}
		\hline
		\multicolumn{2}{|l|}{\# comp.~nodes}                               &                    1 &                    2 &                    4 &                    8 &                   16 \\ \hline\hline
		\multicolumn{2}{|l|}{matrix size}                                  &                 375K &                 750K &                 1.5M &                   3M &                   6M \\ \hline
		\multicolumn{2}{|l|}{CPU}                                          &            {\bf 1.3} &            {\bf 1.6} &            {\bf 1.7} &            {\bf 1.8} & {\bf 1.9} \\ \hline
		\multirow{4}{*}{\rotatebox[origin=c]{90}{GPU}} & $n_p$/gpu = 
   1 &                  3.2 &                  3.5 &                  3.9 &                  5.4 &                  5.6 \\
		                                               & 2                 &                  2.1 &                  2.3 &                  2.9 &                  3.3 &                  3.4 \\
		                                               & 4                 &                  1.4 &                  1.8 &                  2.1 &                  2.2 &                  2.3 \\
		                                               & 6                 &                  1.5 &                  1.6 &                  1.7 &                  2.0 &                  2.3 \\
		                                               & 7                 &            {\bf 1.2} &            {\bf 1.6} &            {\bf 1.7} &            {\bf 2.0} &            {\bf 2.2} \\ \hline
		\multicolumn{2}{|l|}{slowdown}                                     & $\mathbf{0.9\times}$ & $\mathbf{1.0\times}$ & $\mathbf{1.0\times}$ & $\mathbf{1.1\times}$ & $\mathbf{1.1\times}$ \\ \hline
	\end{tabular}
   \caption{Tacho.}
\end{subfigure}
}
\caption{Numerical Setup Time in seconds for Weak-scale 3D elasticity problems on Summit. Number of MPI process per GPU changes between every GPU row from 1 to 7. This improves the numerical setup time up to 17x and 3x for SuperLU and Tacho on GPU runs. The GPU runs are slightly slower than CPU runs in this phase.
         \label{tab:setup}}
\end{table}

\begin{itemize}
\item
{\tt SuperLU} performs partial pivoting during its numerical factorization.
This ensures the numerical stability of the solver, but
the sparsity structures of the LU factors depend on
the numerical values. As a result there is very little work that can be reused
from the symbolic factorization. For instance, with {\tt SuperLU},
both the symbolic and numerical setups for the {\tt Kokkos-Kernels}
sparse-triangular solver need to be performed after each numerical factorization,
which takes up significant part of the difference between the setup times on CPUs and with GPUs in the plot.
\item
On the other hand, {\tt Tacho} performs the pivoting only within its frontal matrices
and given the same sparsity structure of the input matrices,
the sparsity structures of the LU factors stay the same, allowing us to reuse
the symbolic setup for the numerical factorization of different matrices. 
In addition, {\tt Tacho} can utilize the GPU. 
Overall, we see similar numerical setup times of {\tt Tacho}
on CPUs and with GPUs.
\end{itemize}

Table~\ref{tab:setup} compares the the weak-scaling numerical setup time of {\tt FROSch} using up to 672 CPU cores 
and 96 GPUs.
As we discussed in Section~\ref{sec:setups}, the computational cost for factoring
the local overlapping subdomain scales superlinearly to the size of the local matrix.
\emph{Hence, similar to the iterative solution time,
the numerical setup time was reduced significantly using MPS
(obtaining speedups of $15 \sim 17\times$ with {\tt SuperLU} and $2 \sim 3\times$ with {\tt Tacho})}. 
Running multiple MPI processes on each GPU also reduced the memory required to store
the LU factors, enabling the solution of a larger linear system.

\emph{Overall, using {\tt Tacho}, the total solution time (the sum of setup and solve time for solving a single linear system)
was $1.1 \sim 1.8\times$ faster with GPUs.
If the application requires to solve a sequence of linear systems (the same matrix~$A$ but with different
right-hand-side $b$ in sequence), then the cost of the numerical
factorization can be amortized over the multiple solves, and speedups closer to $2\times$ can be obtained.}

\begin{figure}
\centerline{
\begin{subfigure}[b]{.48\linewidth}
  \includegraphics[width=\linewidth]{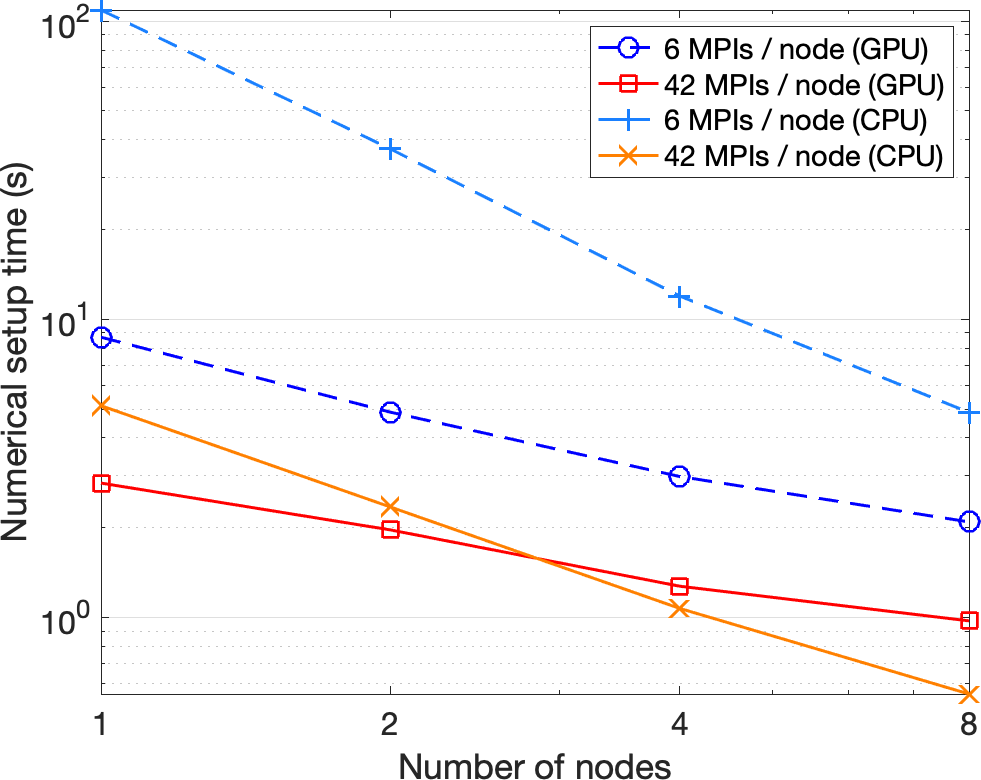}
  \caption{Numerical setup time.}  
\end{subfigure}
\begin{subfigure}[b]{.48\linewidth}
  \includegraphics[width=\linewidth]{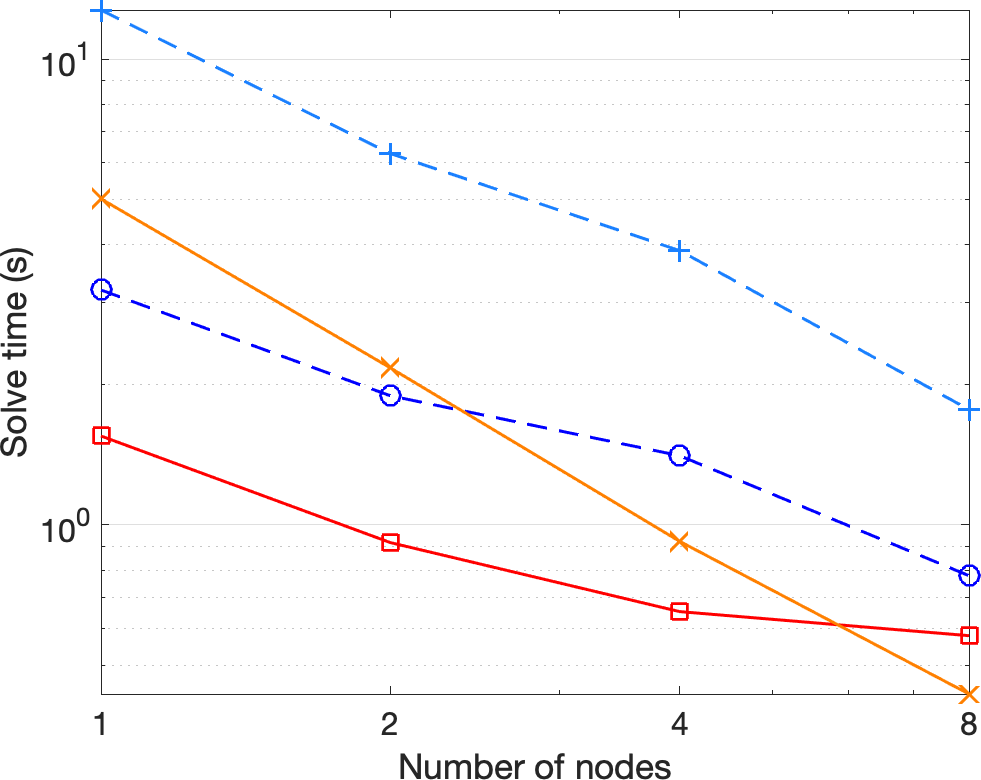}
  \caption{Solve time.}  
\end{subfigure}
}
\caption{Strong parallel scaling with 3D elasticity ($n=1$M).
\label{fig:strong}}
\end{figure}

To summarize our studies with the exact solver, Fig.~\ref{fig:strong} shows the strong parallel-scaling results,
where either 6 or 42 MPI processes were used on each node. For our CPU runs with 6 MPI processes per node,
we linked {\tt Tacho} with the threaded version of ESSL and used 7 threads for each MPI process.
We clearly see the advantage of having 42 MPI processes on each node for both CPU and GPU runs.
\emph{Overall, the GPUs can provide speedups for both setup and solve time
as long as the local matrix sizes are large enough.}

\subsection{Approximate Local Solvers}

\subsubsection{Incomplete LU Factorization}

\begin{table}
\centerline{
\begin{subfigure}[b]{\linewidth} \centering \ssmall \addtolength{\tabcolsep}{-1.2mm}
  \centerline{
   \begin{tabular}{|l|l|rrrr|}
   	\hline
   	\multicolumn{2}{|l|}{ILU level}                           &                    0 &                    1 &                    2 &                    3 \\ \hline\hline
    \multirow{2}{*}{\rotatebox[origin=c]{90}{CPU}} & No            & {\bf 1.5} & {\bf 1.9} & {\bf 3.0} & {\bf 4.8}\\
    & ND            & 1.6 & 2.6 & 4.4 & 7.4\\
   \hline\hline
    \multirow{4}{*}{\rotatebox[origin=c]{90}{GPU}} & KK(No)        & 1.4 & 1.5 & 1.8 & 2.4\\
    & KK(ND)        & 1.7 & 2.0 & 2.9 & 5.2\\
    \cline{2-6}   
    & Fast(No)      & {\bf 1.5} & {\bf 1.6} & {\bf 2.1} & {\bf 3.2}\\ 
    & Fast(ND)      & 1.5 & 1.7 & 2.5 & 4.5\\
    \hline
    \multicolumn{2}{|l|}{speedup}       & $\mathbf{1.0\times}$ & $\mathbf{1.2\times}$ & $\mathbf{1.4\times}$ & $\mathbf{1.5\times}$ \\
    \hline
   \end{tabular}
   \raisebox{-.5\height}{\includegraphics[width=.5\linewidth]{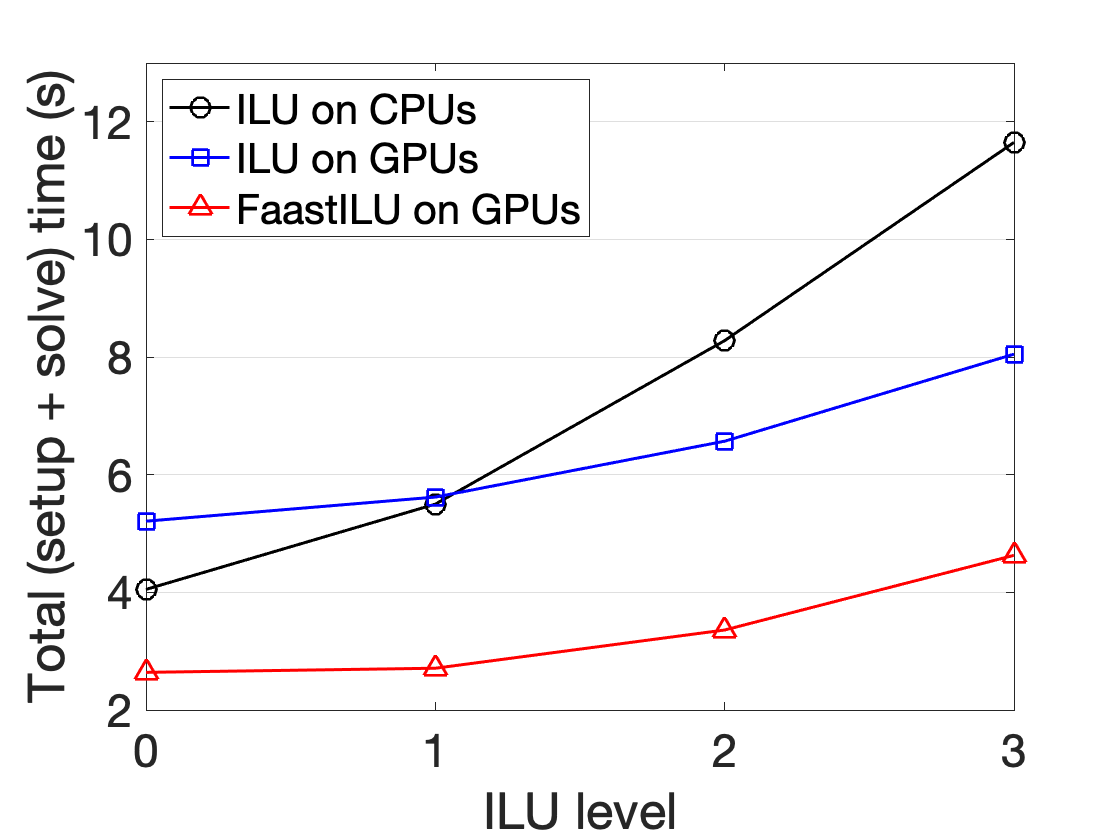}}
   }
   \caption{Setup (table on left) and total time (figure on right).\\\;\smallskip
   }
   \label{tab:ilu-level-setup}
\end{subfigure}
}
\centerline{
\begin{subfigure}[b]{\linewidth} \centering \footnotesize 
   \begin{tabular}{|l|l|rrrr|}
   	\hline
   	\multicolumn{2}{|l|}{ILU level}                           &                    0 &                    1 &                    2 &                    3 \\ \hline\hline
   	\multirow{2}{*}{\rotatebox[origin=c]{90}{CPU}} & No       &     {\bf 2.55 (158)} &     {\bf 3.60 (112)} &  {\bf 5.28\;\; (99)} &  {\bf 6.85\;\; (88)} \\
   	                                               & ND       &           4.17 (227) &           5.36 (134) &           6.61 (105) &        7.68\;\; (88) \\ \hline\hline
   	\multirow{4}{*}{\rotatebox[origin=c]{90}{GPU}} & KK(No)   &           3.81 (158) &           4.12 (112) &        4.77\;\; (99) &        5.65\;\; (88) \\
   	                                               & KK(ND)   &           2.89 (227) &           4.27 (134) &           5.57 (105) &        6.36\;\; (88) \\ \cline{2-6}
   	                                               & Fast(No) &     {\bf 1.14 (173)} &     {\bf 1.11 (141)} &     {\bf 1.26 (134)} &     {\bf 1.43 (126)} \\
   	                                               & Fast(ND) &           1.49 (227) &           1.15 (137) &           1.10 (109) &           1.22 (100) \\ \hline
   	\multicolumn{2}{|l|}{speedup}                             & $\mathbf{2.2\times}$ & $\mathbf{3.2\times}$ & $\mathbf{4.3\times}$ & $\mathbf{4.8\times}$ \\ \hline
   \end{tabular}
   \caption{Solve time  (iteration count).}
   \label{tab:ilu-level-solve}
\end{subfigure}
}
\caption{Performance of {\tt FROSch} for 3D elasticity problems on one Summit node ($n=648$K on 42 MPI processes)
         using local {\tt Kokkos-Kernels} ILU (KK) or {\tt FastILU} (Fast) and
         no reordering (No) or nested dissection (ND).}
\label{tab:ilu-level}
\end{table}

We now study the effects of using an incomplete LU (ILU) factorization as our local solver
on the performance of {\tt FROSch}.
For these experiments, we used \emph{level}-based ILU($k$) as local solver only for solving the local overlapping subdomain problems,
while {\tt Tacho} was used for computing the basis function and for solving the coarse problem.
The inexact solvers reduce the required storage cost, and hence,
we are solving larger linear systems in this section, compared to those solved in Section~\ref{sec:results-exact}.

Table~\ref{tab:ilu-level-setup} shows the effects of the number of ILU levels, $k$, on the numerical setup time.
Besides {\tt SpILU} and {\tt SpTRSV} (based on level-set scheduling),
we also show the performance of their iterative variants, {\tt FastILU} and {\tt FastSpTRSV},
where we performed three and five Jacobi iterations, respectively.
As we increase the level (and the computation required to compute ILU factors increases),
the speedup gained using the GPUs for the numerical setup time increased. 

Table~\ref{tab:ilu-level-solve} shows the total iteration time
with increasing levels for ILU. {\tt FastILU} computes approximation to the ILU factors,
and {\tt FastSpTRSV} solves the triangular system, approximately.
As a result, compared to {\tt SpILU}, 
GMRES required more iterations to converge using {\tt FastILU}.
However, they provide more parallelism, which the GPU can exploit.
\emph{Overall, GMRES had the fastest time to solution using the iterative variants
with the speedups of $2.8\sim4.4\times$.}

\begin{table}
\centerline{
\begin{subfigure}[b]{\linewidth} \centering \footnotesize
   \begin{tabular}{|ll|rrrrr|}
   	\hline
   	\multicolumn{2}{|l|}{\# comp.~nodes}                  &                    1 &                    2 &                    4 &                    8 &                   16 \\ \hline\hline
   	\multicolumn{2}{|l|}{matrix size}                     &                 648K &                 1.2M &                 2.6M &                 5.2M &                10.3M \\ \hline
   	\multicolumn{2}{|l|}{CPU}                             &            {\bf 1.9} &            {\bf 2.2} &            {\bf 2.4} &            {\bf 2.4} &            {\bf 2.6} \\ \hline
   	\multirow{2}{*}{\rotatebox[origin=c]{90}{GPU}} & KK   &                  1.4 &                  2.0 &                  2.2 &                  2.4 &                  2.8 \\
   	                                               & Fast &            {\bf 1.5} &            {\bf 2.2} &            {\bf 2.3} &            {\bf 2.5} &            {\bf 2.8} \\ \hline
   	\multicolumn{2}{|l|}{speedup}                         & $\mathbf{1.3\times}$ & $\mathbf{1.0\times}$ & $\mathbf{1.0\times}$ & $\mathbf{1.0\times}$ & $\mathbf{0.9\times}$ \\ \hline
   \end{tabular}
   \caption{Setup time (s).\\\;\smallskip}
\end{subfigure}
}
\centerline{
\begin{subfigure}[b]{\linewidth} \centering \scriptsize 
   \begin{tabular}{|ll|rrrrr|}
   	\hline
   	\multicolumn{2}{|l|}{\# comp.~nodes}                  &                    1 &                    2 &                    4 &                    8 &                   16 \\ \hline\hline
   	\multicolumn{2}{|l|}{matrix size}                     &                 648K &                 1.2M &                 2.6M &                 5.2M &                10.3M \\ \hline
   	\multicolumn{2}{|l|}{CPU}                             &      {\bf 4.0 (119)} &      {\bf 3.8 (110)} &      {\bf 3.7 (105)} &   {\bf 3.3\;\; (97)} &      {\bf 4.1 (109)} \\ \hline
   	\multirow{2}{*}{\rotatebox[origin=c]{90}{GPU}} & KK   &            4.3 (119) &            3.9 (110) &            4.8 (105) &         4.3\;\; (97) &            4.9 (109) \\
   	                                               & Fast &      {\bf 1.2 (154)} &      {\bf 1.0 (133)} &      {\bf 1.1 (130)} &      {\bf 1.3 (117)} &      {1.6 (131)\bf } \\ \hline
   	\multicolumn{2}{|l|}{speedup}                         & $\mathbf{3.3\times}$ & $\mathbf{3.8\times}$ & $\mathbf{3.4\times}$ & $\mathbf{2.5\times}$ & $\mathbf{2.6\times}$ \\ \hline
   \end{tabular}
   \caption{Solve time (s) (iteration count).}
\end{subfigure}
}
\caption{Weak scaling Parallel Performance (42 MPIs on each node) of 3D elasticity problems on Summit,
         using ILU(1) as local subdomain solvers.}\label{tab:ilu}
\end{table}

Finally,
Table~\ref{tab:ilu} shows weak-scaling results using the inexact ILU(1) local solver on up to 672 CPU cores 
and 96 GPUs.
For all these experiments, we used the original matrix ordering since the matrix reordering did not improve the performance significantly, while it could increase the iteration count.
Even with the inexact local solver,
the iteration counts were 
almost independent
of the number of subdomains. It can be seen in Table~\ref{tab:ilu} that even with the higher iteration count, the inexact (Fast) option is faster than the exact triangular solve (KK). \emph{We see $3.1 \sim 4.4\times$ speedups using the iterative variants on GPUs. We also observe $2.5 \sim3.8\times$ speedup using GPUs. The setup times are nearly the same on CPUs and GPUs with ILU(1) on multiple nodes.}

\subsubsection{Single-precision FROSch}

\begin{table}
\centerline{
\begin{subfigure}[b]{\linewidth} \centering \footnotesize 
   \begin{tabular}{|ll|rrrrr|}
   	\hline
   	\multicolumn{2}{|l|}{\# comp.~nodes}                     &           1 &           2 &           4 &           8 &          16 \\ \hline\hline
   	\multicolumn{2}{|l|}{matrix size}                        &        375K &        750K &        1.5M &          3M &          6M \\ \hline
   	\multirow{3}{*}{\rotatebox[origin=c]{90}{CPU}} & double  &         2.5 &         3.0 &         3.3 &         3.8 &         3.6 \\
   	                                               & single  &         1.8 &         2.2 &         2.3 &         2.6 &         2.6 \\
   	                                               & speedup & $1.4\times$ & $1.4\times$ & $1.4\times$ & $1.5\times$ & $1.4\times$ \\ \hline
   	\multirow{3}{*}{\rotatebox[origin=c]{90}{GPU}} & double  &         3.5 &         4.2 &         4.8 &         5.6 &         5.4 \\
   	                                               & single  &         2.6 &         3.2 &         3.4 &         4.0 &         4.0 \\
   	                                               & speedup & $1.3\times$ & $1.3\times$ & $1.4\times$ & $1.4\times$ & $1.4\times$ \\ \hline
   \end{tabular}
   \caption{SuperLU.\\\;\smallskip}
\end{subfigure}
}
\centerline{
\begin{subfigure}[b]{\linewidth} \centering \footnotesize 
   \begin{tabular}{|ll|rrrrr|}
   	\hline
   	\multicolumn{2}{|l|}{\# comp.~nodes}                     &           1 &           2 &           4 &           8 &          16 \\ \hline\hline
   	\multicolumn{2}{|l|}{matrix size}                        &        375K &        750K &        1.5M &          3M &          6M \\ \hline
   	\multirow{3}{*}{\rotatebox[origin=c]{90}{CPU}} & double  &         1.3 &         1.6 &         1.7 &         1.8 &         1.9 \\
   	                                               & single  &         1.0 &         1.2 &         1.3 &         1.4 &         1.4 \\
   	                                               & speedup & $1.3\times$ & $1.3\times$ & $1.3\times$ & $1.3\times$ & $1.4\times$ \\ \hline
   	\multirow{3}{*}{\rotatebox[origin=c]{90}{GPU}} & double  &         1.2 &         1.6 &         1.7 &         2.0 &         2.2 \\
   	                                               & single  &         1.0 &         1.3 &         1.4 &         1.7 &         2.0 \\
   	                                               & speedup & $1.2\times$ & $1.2\times$ & $1.2\times$ & $1.2\times$ & $1.1\times$ \\ \hline
   \end{tabular}
   \caption{Tacho.}
\end{subfigure}
}
\caption{Numerical Setup Time in seconds, using single or double precision {\tt FROSch}: for Weak-scale 3D elasticity problems on Summit.}
\end{table}

\begin{table}
\centerline{
\begin{subfigure}[b]{\linewidth} \centering \footnotesize \addtolength{\tabcolsep}{-0.5mm}
   \begin{tabular}{|ll|rrrrr|}
   	\hline
   	\multicolumn{2}{|l|}{\# comp.~nodes}                   &           1 &           2 &           4 &           8 &          16 \\ \hline\hline
   	\multicolumn{2}{|l|}{matrix size}                        &        375K &        750K &        1.5M &          3M &          6M \\ \hline
   	\multirow{3}{*}{\rotatebox[origin=c]{90}{CPU}} & double  &   2.03 (75) &   2.07 (69) &   1.87 (61) &   1.95 (58) &   2.48 (69) \\
   	                                               & single  &   1.89 (76) &   1.60 (69) &   1.71 (62) &   1.75 (58) &   2.37 (69) \\
   	                                               & speedup & $1.0\times$ & $1.3\times$ & $1.1\times$ & $1.1\times$ & $1.0\times$ \\ \hline
   	\multirow{3}{*}{\rotatebox[origin=c]{90}{GPU}} & double  &   1.03 (75) &   1.04 (69) &   0.90 (61) &   0.97 (58) &   1.18 (69) \\
   	                                               & single  &   1.01 (75) &   1.03 (69) &   0.98 (62) &   1.10 (58) &   1.28 (69) \\
   	                                               & speedup & $1.0\times$ & $1.0\times$ & $0.9\times$ & $0.9\times$ & $0.9\times$ \\ \hline
   \end{tabular}
   \caption{SuperLU.\\\;\smallskip}
\end{subfigure}
}
\centerline{
\begin{subfigure}[b]{\linewidth} \centering \footnotesize \addtolength{\tabcolsep}{-0.5mm}
   \begin{tabular}{|ll|rrrrr|}
   	\hline
   	\multicolumn{2}{|l|}{\# comp.~nodes}                   &           1 &           2 &           4 &           8 &          16 \\ \hline\hline
   	\multicolumn{2}{|l|}{matrix size}                        &        375K &        750K &        1.5M &          3M &          6M \\ \hline
   	\multirow{3}{*}{\rotatebox[origin=c]{90}{CPU}} & double  &   1.60 (75) &   1.63 (69) &   1.49 (61) &   1.51 (58) &   1.90 (69) \\
   	                                               & single  &   1.11 (76) &   1.13 (69) &   1.02 (62) &   1.04 (58) &   1.30 (69) \\
   	                                               & speedup & $1.4\times$ & $1.4\times$ & $1.4\times$ & $1.4\times$ & $1.5\times$ \\ \hline
   	\multirow{3}{*}{\rotatebox[origin=c]{90}{GPU}} & double  &   0.99 (75) &   0.93 (69) &   0.82 (61) &   0.93 (58) &   1.22 (69) \\
   	                                               & single  &   1.00 (75) &   0.92 (69) &   0.84 (62) &   0.93 (58) &   1.21 (69) \\
   	                                               & speedup & $1.0\times$ & $1.0\times$ & $1.0\times$ & $1.0\times$ & $1.0\times$ \\ \hline
   \end{tabular}
   \caption{Tacho.}
\end{subfigure}
}
\caption{Total iteration time and Iteration count, using double or single precision {\tt FROSch} for Weak-scale 3D elasticity problems on Summit.
\label{tab:single}}
\end{table}

Even though typical scientific applications require double precision accuracy,
some emerging hardware delivers lower-precision arithmetic at higher performance.
There are other machines that provide the same performance for double and single
precision arithmetic.
However, even in that case, using a lower-precision arithmetic
reduces the required amount of data transfer.
Since the performance of the sparse solver is often limited by the memory bandwidth,
reducing the required communication volume alone could reduce the solver time.

Table~\ref{tab:single} shows the performance results,
where {\tt FROSch} in single precision is used to precondition GMRES in double precision.
For these particular problems, the setup time was reduced using single-precision {\tt FROSch}, 
while the number of iterations required for the convergence \emph{to same accuracy as double precision use cases} is maintained. \emph{Specifically, using single precision, both on 672 CPU cores and 96 GPUs, we observe $1.3 \sim 1.5\times$ speedup in {\tt SuperLU} based setup, while we observe $1.1 \sim 1.4\times$ speedup in {\tt Tacho} based setup.} We do not see a benefit in solve times when using single precision (Table~\ref{tab:single}). \emph{Nevertheless, GMRES converges in similar number of iterations when using single or double precision preconditioner.}

\subsection{Summary of Key Results}

Using multiple subdomains per GPU improves the performance considerably. In terms of the setup time using a direct factorization for the subdomain solver, GPU-based factorization in {\tt Tacho} provides a distinct advantage over CPU-based {\tt SuperLU}. There is no distinct advantage in the solve time on the GPUs when using {\tt SuperLU} and {\tt Kokkos Kernels} combination or using {\tt Tacho}. Incomplete factorizations allow us to solve larger problems. Though we do not see noticeable difference in performance over direct factorization due to the trade-off between number of iterations and setup/solve times, the solve time was reduced using GPUs. Iterative incomplete factorizations and triangular solve result in significant speedup compared to the standard incomplete factorizations even with increased number of iterations.  Using lower precision computations allows us to solve a larger linear system, and improved the setup time though not the solve time.

%% file: sections/concl.tex
We presented {\tt FROSch}, which
implements the GDSW algorithm for GPU cluster within {\tt Trilinos} software framework.
Our performance results on Summit supercomputer with NVIDIA V100 GPUs
demonstrated the potential of {\tt FROSch}:
with GPUs, the numerical setup times remain about the same as that on CPUs, while the solve time
can be reduced by factors of around $2\times$. We presented a thorough experimental study varying several solver options from two direct solvers, incomplete factorization techniques with different level of fill and different orderings, inexact factorizations and the use of lower precision arithmetic. 
Though we only showed the performance results with NVIDIA GPUs,
our implementation is portable to other GPUs through the use of Kokkos.
We plan to study performance of {\tt FROSch}
with AMD GPUs.